\documentclass[11pt]{article}
\usepackage{amsfonts,amsmath,amsxtra}
\usepackage{latexsym}
\usepackage[matrix,arrow,curve]{xy}
\usepackage{amssymb}

\def\hybrid{\topmargin 0pt      \oddsidemargin 0pt
        \headheight 0pt \headsep 0pt
        \textwidth 16.5cm
        \textheight 23.5cm
        \marginparwidth 0.0in
        \parskip 5pt plus 1pt   \jot = 1.5ex}
\catcode`\@=11
\def\marginnote#1{}
\newcount\hour
\newcount\minute
\newtoks\amorpm
\hour=\time\divide\hour by60 \minute=\time{\multiply\hour by60
\global\advance\minute by-\hour}
\edef\standardtime{{\ifnum\hour<12 \global\amorpm={am}%
        \else\global\amorpm={pm}\advance\hour by-12 \fi
        \ifnum\hour=0 \hour=12 \fi
      \number\hour:\ifnum\minute<10 0\fi\number\minute\the\amorpm}}
\edef\militarytime{\number\hour:\ifnum\minute<10 0\fi\number\minute}

\def\draftlabel#1{{\@bsphack\if@filesw {\let\thepage\relax
   \xdef\@gtempa{\write\@auxout{\string
      \newlabel{#1}{{\@currentlabel}{\thepage}}}}}\@gtempa
   \if@nobreak \ifvmode\nobreak\fi\fi\fi\@esphack}
        \gdef\@eqnlabel{#1}}
\def\@eqnlabel{}
\def\@vacuum{}
\def\draftmarginnote#1{\marginpar{\raggedright\scriptsize\tt#1}}

\def\draft{\oddsidemargin -0.1truein
        \def\@oddfoot{\sl preliminary draft \hfil
        \rm\thepage\hfil\sl\today\quad\militarytime}
        \let\@evenfoot\@oddfoot \overfullrule 3pt
        \let\label=\draftlabel
        \let\marginnote=\draftmarginnote
\def\@eqnnum{{\rm (\theequation)}
\rlap{\kern\marginparsep\tt\@eqnlabel}%
\global\let\@eqnlabel\@vacuum}  }

\newfont{\Bbbb}{msbm7 scaled 1\@ptsize00}
\newcommand{\zs}{\raise-1pt\hbox{$\mbox{\Bbbb Z}$}}

scaled 1\@ptsize00

\font\sevenmsa=msam6 
\newfam\msafam
\textfont\msafam=\sevenmsa
\def\hexnumber@#1{\ifnum#1<10 \number#1\else
\ifnum#1=10 A\else\ifnum#1=11 B\else\ifnum#1=12 C\else \ifnum#1=13
D\else\ifnum#1=14 E\else\ifnum#1=15 F\fi\fi\fi\fi\fi\fi\fi}
\def\msa@{\hexnumber@\msafam}
\def\llcorner{\delimiter"4\msa@78\msa@78 }
\def\lrcorner{\delimiter"5\msa@79\msa@79 }
\mathchardef\blacktriangleright="3\msa@49
\mathchardef\blacktriangleleft="3\msa@4A \font\tenmsb=msbm10 scaled
1\@ptsize00
\newfam\msbfam
\textfont\msbfam=\tenmsb \scriptfont\msbfam=\tenmsb

\newdimen\Squaresize \Squaresize=14pt
\newdimen\Thickness \Thickness=0.5pt

\def\Square#1{\hbox{\vrule width \Thickness
   \vbox to \Squaresize{\hrule height \Thickness\vss
      \hbox to \Squaresize{\hss#1\hss}
   \vss\hrule height\Thickness}
\unskip\vrule width \Thickness} \kern-\Thickness}

\def\Vsquare#1{\vbox{\Square{$#1$}}\kern-\Thickness}

\def\numberbysection{\@addtoreset{equation}{section}
        \def\theequation{\thesection.\arabic{equation}}}
\numberbysection

\renewcommand{\theequation}{\thesection.\arabic{equation}}
\def\titlepage{\@restonecolfalse\if@twocolumn\@restonecoltrue\onecolumn
     \else \newpage \fi \thispagestyle{empty}\c@page\z@
        \def\thefootnote{\fnsymbol{footnote}} }

\def\endtitlepage{\if@restonecol\twocolumn \else  \fi
        \def\thefootnote{\arabic{footnote}}
        \setcounter{footnote}{0}}  
\relax

\hybrid
\makeatletter
\newdimen\normalarrayskip            
\newdimen\minarrayskip               
\normalarrayskip\baselineskip \minarrayskip\jot
\newif\ifold             \oldtrue            \def\new{\oldfalse}
\def\arraymode{\ifold\relax\else\displaystyle\fi}
\def\eqnumphantom{\phantom{(\theequation)}} 
\def\@arrayskip{\ifold\baselineskip\z@\lineskip\z@
     \else
     \baselineskip\minarrayskip\lineskip1\baselineskip\fi}

\def\@arrayclassz{\ifcase \@lastchclass \@acolampacol \or
\@ampacol \or \or \or \@addamp \or
   \@acolampacol \or \@firstampfalse \@acol \fi
\edef\@preamble{\@preamble
  \ifcase \@chnum
     \hfil$\relax\arraymode\@sharp$\hfil
     \or $\relax\arraymode\@sharp$\hfil
     \or \hfil$\relax\arraymode\@sharp$\fi}}

\def\@array[#1]#2{\setbox\@arstrutbox=\hbox{\vrule
     height\arraystretch \ht\strutbox
     depth\arraystretch \dp\strutbox
width\z@}\@mkpream{#2}\edef\@preamble{\halign \noexpand\@halignto
\bgroup \tabskip\z@ \@arstrut \@preamble \tabskip\z@ \cr}%
\let\@startpbox\@@startpbox \let\@endpbox\@@endpbox
    \if #1t\vtop \else \if#1b\vbox \else \vcenter \fi\fi
  \bgroup \let\par\relax
  \let\@sharp##\let\protect\relax
  \@arrayskip\@preamble}
\def\eqnarray{\stepcounter{equation}%
              \let\@currentlabel=\theequation
              \global\@eqnswtrue
              \global\@eqcnt\z@
              \tabskip\@centering              
              \let\\=\@eqncr
              $$%
            \halign to \displaywidth  \bgroup
             \eqnumphantom \@eqnsel
      \hskip\@centering                               
    $\displaystyle  \tabskip\z@ {##}$%
    &\global\@eqcnt\@ne \hskip 2\arraycolsep
         $ \displaystyle  \arraymode{##}$\hfil
    &\global\@eqcnt\tw@ \hskip 2\arraycolsep
         $\displaystyle\tabskip\z@{##}$\hfil
         \tabskip\@centering
    &{##}\tabskip\z@\cr}
\makeatother

\newcommand{\CC}{{\mathbb{C}}}

\def\IC{\mathbb{C}}

\def\IP{\mathbb{P}}

\def\IR{\mathbb{R}}
\def\IZ{\mathbb{Z}}

\def\CC {\mathcal{C}}
\def\CD {\mathcal{D}}

\def\CF {\mathcal{F}}

\def\CH {\mathcal{H}}
\def\CI {\mathcal{I}}

\def\CL {\mathcal{L}}
\def\CM {\mathcal{M}}

\def\CT {\mathcal{T}}
\def\CU {\mathcal{U}}
\def\CV {\mathcal{V}}
\def\CW {\mathcal{W}}

\def\CZ {\mathcal{Z}}

\def\frak{\mathfrak}

\def\we{\raise-1pt\hbox{$\,\stackrel{\wedge}{,}\,$}}
\def\tr{{\rm tr}\,}
\def\pr {\partial}

\def\ov {{\overline}}
\def\a {{\alpha}}
\def\be {{\beta}}

\def\s {{\sigma}}

\def\la{\lambda}

\def\e{\epsilon}
\def\pr {\partial}

\newcommand{\ul}{{\underline{\la}}}

\newcommand{\At}{{{}^\top\!A}}
\newcommand{\Bt}{{{}^\top\!B}}
\newcommand{\Dt}{{{}^\top\!D}}
\newcommand{\Pt}{{{}^\top\!P}}

\newcommand{\upp}{{\underline{p}}}

\newcommand{\ux}{{\underline{x}}}

\newcommand{\gl}{{\mathfrak{gl}}}

\newcommand{\Fl}{{\mathrm{ Fl}}}
\newcommand{\Gr}{{\mathrm{ Gr}}}
\newcommand{\h}{{\hbar}}

\newcommand{\Lie}{{\mathrm{ Lie\ }}}

\newcommand{\Mat}{{\mathrm{ Mat}}}

\newtheorem{te}{Theorem}[section]
\newtheorem{de}{Definition}[section]
\newtheorem{prop}{Proposition}[section]
\newtheorem{cor}{Corollary}[section]
\newtheorem{lem}{Lemma}[section]

\newtheorem{rem}{Remark}[section]

\newcommand\bqa{\begin{eqnarray}}
\newcommand\eqa{\end{eqnarray}}
\def\be{\begin{eqnarray}\new\begin{array}{cc}}
\def\ee{\end{array}\end{eqnarray}}

\def\beq{\begin{equation}}
\def\eeq{\end{equation}}
\def\bse{\begin{subequations}}                
\def\ese{\end{subequations}}
\def\bp{\begin{pmatrix}}
\def\ep{\end{pmatrix}}

\def\h{\hbar}
\def\i{\imath}

\newcommand\rk{\operatorname{rank}}

\def\stack#1#2{\raise0.7pt\hbox{$\mathrel{\mathop{#2}\limits^{#1}}$}}
\def\tr{\triangleright}
\def\tl{\triangleleft}
\def\sem{\mathsurround=0pt \raise1pt
\hbox{$\scriptscriptstyle>\!\!$}\:\!\!\tl}
\def\mes{\mathsurround=0pt \tr\!\:\!\raise0.8pt
\hbox{$\scriptscriptstyle\!\!<$}\,}
\def\]{\mathsurround=0pt ]\raise-2pt\hbox{$_\ast$}}

\def\<{\langle}
\def\>{\rangle}
\def\ov{\overline}

\newcounter{pac}[section]
\newcommand{\npa}{\addtocounter{pac}{1} \noindent {\bf
\arabic{section}.\arabic{pac}}\,\,\,}

\newcounter{pacc}[subsection]
\newcommand{\npaa}{\addtocounter{pacc}{1} \noindent {\bf
\arabic{section}.\arabic{subsection}.\arabic{pacc}}\,\,\,}
  0

\begin{document}

\title{\bf On parabolic Whittaker functions II}
\author{\sc\Large Sergey Oblezin\footnote{The work is
partially supported by P. Deligne's 2004 Balzan Prize in
Mathematics, and by the RFBR 09-01-93108-NCNIL-a Grant.}}

\date{}
\maketitle

\begin{abstract}\noindent
We derive a Givental-type stationary phase integral representation
for the specified $\Gr_{m,N}$-Whittaker function introduced in
\cite{GLO2}, which presumably describes the $S^1\times
U_N$-equivariant Gromov-Witten invariants of Grassmann variety
$\Gr_{m,N}$. Our main tool is a generalization of Whittaker model
for principal series $\CU(\frak{gl}_N)$-modules. In particular, our
construction includes a representation theory interpretation of the
Batyrev--Ciocan-Fontanine--Kim--van Straten toric degeneration of
Grassmannian, providing a direct connection between this toric
degeneration of $\Gr_{m,N}$ and total positivity for unipotent
matrices.
\end{abstract}

\section*{Introduction}

Quantum cohomology D-module $QH^*(\Fl_N)$ for complete flag manifold
$\Fl_N$ of $GL_N$ can be identified with the quantum Toda lattice
\cite{Giv1}, \cite{GK}. In \cite{Giv2} Givental proposed a
stationary phase integral formula for the solution (generating
function) of $QH^*(\Fl_N)$. Unfortunately, Givental's approach
cannot be applied directly to description of quantum cohomology of
incomplete (partial) flag variaties, since no relevant Toda lattice
associated with a partial flag variaties was known.

Givental's construction possesses a set of remarkable properties,
and particularly it involves a degeneration of $\Fl_N$ to a certain
Gorenstein toric Fano variety. In \cite{BCFKvS1} and \cite{BCFKvS2}
Batyrev, Ciocan-Fontanine, Kim and van Straten introduced toric
degenerations of partial flag manifolds of $GL_N$, generalizing
Givental's result; they conjectured that the proposed toric
varieties describe the quantum cohomology of partial flag manifolds.

In \cite{GKLO}, \cite{GLO1} was given a representation theory proof
of Givental's stationary phase integral formula for the complete
flag manifold, including a representation theory derivation of
Givental's toric degeneration of complete flag manifold in more
general context of Lie groups of classical type. In this note we
extend the construction of \cite{GKLO} to Grassmann varieties,
giving a representation theory construction of the toric
degeneration of Grassmann varieties $\Gr_{m,N}$ proposed by Batyrev,
Ciocan-Fontanine, Kim, van Straten in \cite{BCFKvS1} and
\cite{BCFKvS2}.

The key obstacle in representation theory approach to quantum
cohomology of homogeneous spaces until recently was an absence of a
relevant Toda lattice associated with partial flag manifolds.
Recently such a quantum Toda-type lattice was proposed in
\cite{GLO2}, using a generalization of the classical Whittaker model
for principal series representations of $\CU(\gl_N)$. The solution
$\Psi^{(m,N)}_{\ul}(x_1,\ldots,x_N)$ to the quantum parabolic Toda
lattice was defined in \cite{GLO2} as a certain matrix element in
principal series representation, and it was referred to as parabolic
Whittaker function, or $\Gr_{m,N}$-Whittaker function. It was
conjectured in \cite{GLO2} that the specialized parabolic Whittaker
function, $\Psi^{(m,N)}_{\ul}(x,0,\ldots,0)$, describes the
$S^1\times U_N$-equivariant quantum cohomology of $\Gr_{m,N}$, and
this conjecture was verified in the case of projective variety
$\Gr_{1,N}=\IP^{N-1}$. This conjecture is supported by an
observation that the specialized symbol
$L(x,0,\ldots,0;\,p_1,\ldots,p_N)$ of the Lax operator associated
with the parabolic Toda lattice reproduces the small cohomology
algebra $qH^*(\Gr_{m,N})$ in the form of \cite{AS} and \cite{K}.

In our main result, Theorem 1.1, we propose a stationary phase
integral formula for the specialized $\Gr_{m,N}$-Whittaker function
$\Psi^{(m,N)}_{\ul}(x,0,\ldots,0)$. Our basic tool is a
generalization of the classical Whittaker model for principal series
representations of $\CU(\gl_N)$, which extends the construction of
generalized Whittaker model from \cite{GLO2}. This result should be
compared with the Mellin-Barnes integral representation for
$\Psi^{(m,N)}_{\ul}(x,0,\ldots,0)$ constructed recently in \cite{O}.
The coincidence of the two integrals in the case
$\Gr_{1,N}=\IP^{N-1}$ is evident due to a simple transformation of
integration variables. For generic $\Gr_{m,N}$ an identification of
the two integrals involves a more delicate description of the
integration contour of the stationary phase integral than we do in
this note; this will be done in a separate note.

Moreover, according to \cite{GLO2} the two integral representations
for $\Psi^{(m,N)}_{\ul}(x,0,\ldots,0)$, the stationary phase one and
the Mellin-Barnes one from \cite{O}, should be presumably identified
with certain correlation functions in the two, type B and type A
respectively, mirror dual topological field theories. A detailed
analysis of the two integral representations of the
$\Gr_{m,N}$-Whittaker function in the framework of \cite{GLO2} will
be given elsewhere.

Our stationary phase integral representation of
$\Gr_{m,N}$-Whittaker function is naturally related to a toric
degeneration of the Grassmannian (Corollary 1.2). Namely, we
identify the function defining the toric variety \footnote{In
\cite{BCFKvS1} it referred to as "Lax operator of Grassmannian"; one
should not mistake this function with the matrix-valued differential
operator $\CL(x,\pr_x)$ in the present work} (in \cite{GLO2} it is
given by the superpotential in the mirror dual Landau-Ginzburg
model) with phase function $\CF_{m,N}$ in our integral formula when
the $U_N$-equivariance parameters $\ul=(\la_1,\ldots,\la_N)$ are
specialized to zero.

Our construction of stationary phase integral uses the
Gauss-Givental realization of universal enveloping algebra
$\CU(\gl_N)$ in the space of functions on totally positive unipotent
(upper-triangular) matrices introduced in \cite{GKLO}. It should be
emphasized that in our representation theory approach the involved
toric degeneration of the (partial) flag manifolds is directly
provided by the total positivity property of the Gauss-Givenatal
realization of $\CU(\gl_N)$; this connection was established in
\cite{GLO1} for complete flag manifolds of Lie groups of classical
types. One may also refer \cite{GLO1} for a detailed analysis of
connections between Gauss-Givental realizations of $\CU(\frak{g})$
in general setting and total positivity phenomenon \cite{L}.

The paper is organized in the following way. In Section 1 we
introduce a generalization of the standard Whittaker model, and
define $\Gr_{m,N}$-Whittaker function.In Theorem 1.1 we propose a
stationary phase integral representation for
$\Psi_{\ul}^{(m,N)}(x,0,\ldots,0)$, which extends the Givental's
integral formula for complete flag manifold $\Fl_N$. Next we
describe the arising toric degeneration of $\Gr_{m,N}$ in Corollary
1.2, and identify it with the Batyrev--Ciocan-Fontanine--Kim--van
Straten construction \cite{BCFKvS1}, \cite{BCFKvS2}. In Proposition
1.1 we derive the Lax operator for the associated $\Gr_{m,\,N}$-Toda
lattice, and in Proposition 1.2 we find out the first two quantum
Hamiltonians of the $\Gr_{m,\,N}$-Toda lattice. The second half
contains detailed proofs of our main results. In particular, in
Section 2 we review on the Gauss-Givental realization of
$\CU(\gl_N)$ from \cite{GKLO} and construct the
$\Gr_{m,N}$-Whittaker vectors; and in Section 3 we verify the
integral formula from Theorem 1.1. In Section 4 we prove
Propositions 1.1 and 1.2.

{\em Acknowledgments}: The author is thankful to A. Gerasimov and D.
Lebedev for very useful discussions.

\newpage

\section{The $\Gr_{m,N}$-Whittaker functions}

A Whittaker model for principal series representation of universal
enveloping alegbra $\CU(\frak{g})$ involves a choice of an
$N$-dimensional commutative subalgebra $\frak{a}$ and the Lie
algebra decomposition
$\frak{g}=\frak{u}_-\oplus\frak{a}\oplus\frak{u}_+$. Then choosing a
pair of characters $\chi_{\pm}:\frak{u}_{\pm}\to\IC$ one can reduce
the space of functions on $G=\Lie(\frak{g})$ to the space of
functions on the commutative subgroup $\Lie(\frak{a})$. In the case
$\frak{g}=\frak{gl}_N$ one uses the Gauss triangular decomposition
of $\frak{gl}_N$ with $\frak{a}$ being the diagonal subalgebra of
semisimple elements, and $\frak{u}_+$ being the nilpotent
subalgebras of upper-triangular matrices.\footnote{Alternatively,
one can consider the Iwasawa decomposition
$\gl_N=\frak{k}\oplus\frak{a}\oplus\frak{u}_+$ with
$\frak{k}\subset\gl_N$ being the compact subalgebra; then choosing
characters of $\frak{k}$ and $\frak{u}_+$ one can obtain another
reduction of the space of functions on $GL_N$}

In \cite{GLO2} was proposed a generalization of the standard scheme
with $\frak{a}\subset\gl_N$ being a commutative subalgebra,
containing both semisimple and nilpotent elements, and similarly for
$\frak{u}_{\pm}$. In the following we use another generalized
Whittaker model, using another choice of commutative subalgebra
$\frak{a}=\frak{h}^{(m,N)}\subset\gl_N$. In special case of
$\Gr_{1,N}=\IP^{N-1}$ our choice of the commutative subalgebra
$\frak{h}^{(m)}\subset\gl_N$ coincides with the one used in
\cite{GLO2} and \cite{O}. For generic $m$ the two commutative
subalgebras differ by a simple automorphism, and thus the
generalized Whittaker model used in our calculations below is
different from \cite{GLO2}.

Let us consider the triangular decomposition of the Lie algebra
$\gl_N$ of real $(N\times N)$-matrices with the standard set of
generators $E_{ij},\,i,j=1,\ldots,N$. Namely, let $\frak{b}_{\pm}$
be the subalgebras of upper- and lower-triangular matrices, and
$\frak{n}_{\pm}=[\frak{b}_{\pm},\,\frak{b}_{\pm}]\subset\frak{b}_{\pm}$
are the radicals of the triangular subalgebras. Then the following
holds:
$$
 \gl_N\,=\,\frak{n}_-\oplus\frak{h}\oplus\frak{n}_+\,.
$$

Next, let $\Delta=\{\alpha_i;\,i\in I\}$ be the set of simple roots
indexed by $I=\{1,2,\ldots,N-1\}$, and $R_+$ be the set of positive
roots. Given an elementary subset $\{m\}\subset I$, let us associate
with $m$ the following modification of the triangular decomposition:
 \be\label{ModTriDecomposition}
  \gl_N\,=\,\frak{n}^{(m,N)}_-\oplus\frak{h}^{(m,N)}\oplus\frak{n}^{(m,N)}_+\,,
 \ee
where the commutative subalgebra $\frak{h}^{(m,N)}$ is spanned by
 \be\label{ModTriDecompositionH}
  H_1\,=\,E_{11}+\ldots+E_{mm},
 \hspace{2.5cm}
  H_i\,=\,E_{i,1}\,,\quad i=2,\ldots,m\,;\\
  H_j\,=\,E_{j,\,N}\,,\quad j=m+1,\ldots,N-1\,;
 \hspace{1.5cm}
  H_N\,=\,E_{m+1,\,m+1}+\ldots+E_{NN}\,.
 \ee
The Lie subalgebras $\frak{n}^{(m)}_{\pm}$ are defined by their set
of generators as follows:
 \bqa\label{ModTriDecompositionN-}
  \frak{n}^{(m,N)}_-\,=\,\Big\<\,E_{m+1,\,1}\,;\quad &
  E_{ki}\,,\,\,i=2,\ldots,m,\,\,k=i,\ldots,N\,; &\\
 \nonumber
 & \hspace{1cm} E_{j+1,\,j}\,,\,\,j=m+1,\ldots,N-1\,\Big\>\,, &
 \eqa
and
 \bqa\label{ModTriDecompositionN+}
  \frak{n}^{(m,N)}_+\,
  =\,\Big\<\,E_{i-1,\,i}\,,\,\,i=2,\ldots,m\,; & \quad
  E_{m,\,N}\,; &\\
 \nonumber
 & E_{kj}\,,\,\,j=m+1,\ldots,N-1,\,\,k=1,\ldots,j\,
 \Big\>\,. &
 \eqa
One may note that
$$
 \dim\frak{h}^{(m,N)}\,=\,\rk\,\frak{gl}_N=N\,,
\hspace{1.5cm}
 \dim\frak{n}^{(m,N)}_{\pm}\,=\,\dim\frak{n}_{\pm}\,
 =\,\frac{N(N-1)}{2}\,.
$$
Let $H^{(m,N)}$ and $N_+^{(m,N)}$ be the  Lie groups corresponding
to the Lie algebras $\mathfrak{h}^{(m,N)}$ and
$\mathfrak{n}_+^{(m,N)}$, then an open part $GL_N^{\circ}$ (the big
Bruhat cell) of $GL_N$ allows the following modification of the
Gauss decomposition:
 \be\label{modGauss}
  GL_N^{\circ}\,\,=\,\,N_-^{(m,N)}\cdot H^{(m,N)}\cdot N_+^{(m,N)}\,.
 \ee

A principal series representation $\CV_{\ul}$ admits a natural
structure of $\CU$-module. Let us assume that the action of the
commutative subalgebra $\frak{h}^{(m)}\subset\gl_N$ in $\CV_{\ul}$
is integrable to the action of commutative subgroup
$H^{(m,N)}\subset GL_N(\IR)$. Below we introduce a pair of elements,
$\psi_L,\,\psi_R\in\CV_{\ul}$, generating a pair of Whittaker
submodules in $\CV_{\ul}$, $\CW_L=\CU\psi_L$ and $\CW_R=\CU\psi_R$.
\begin{de}
The $\Gr_{m,N}$-Whittaker vectors $\psi_L$ and $\psi_R$ are defined
by the following conditions:
 \be\label{LeftWhittEqs}
 \left\{
 \begin{array}{lcc}
  E_{m+1,\,1}\psi_L\,=\,\h^{-1}\psi_L\,; & &\\
  E_{ki}\psi_L\,=\,0\,, & & i=2,\ldots,m,\,\,k=i,\ldots,N\,;\\
  E_{j+1,\,j}\psi_L\,=\,\h^{-1}\psi_L\,, & &
  j=m+1,\ldots,N-1\,;
 \end{array}
 \right.
 \ee
 \be\label{RightWhittEqs}
 \left\{
 \begin{array}{lcc}
  E_{i-1,\,i}\psi_R\,=\,-\h^{-1}\psi_R\,, & & i=2,\ldots,m\,;\\
  E_{kj}\psi_R\,=\,0\,,& & j=m+1,\ldots,N-1,\,\,k=1,\ldots,j\,;\\
  E_{kj}\psi_L\,=\,(-1)^{\e(m,N)}\h^{-1}\psi_L\,, & &
  j=m+1,\ldots,N-1,\,\,k=1,\ldots,j\,;
 \end{array}
 \right.
 \ee
where  $\epsilon(m,N)$ is an integer number and $\h$ is an
indeterminant.
\end{de}
\begin{lem}
The introduced $\Gr_{m,N}$-Whittaker vectors define characters of
the subalgebras $\frak{n}^{(m,N)}_{\pm}$:
 \be
  \chi_+\,:\quad\frak{n}^{(m,N)}_+\,\longrightarrow\,\IC\,,
 \hspace{3.5cm}
  \chi_-\,:\quad\frak{n}^{(m,N)}_-\,\longrightarrow\,\IC\,,
 \ee
\end{lem}
\emph{Proof:} One can readily check that the defining equations
\eqref{LeftWhittEqs} are compatible with Lie algebra relations in
$\frak{n}^{(m,N)}_-$; and the same is valid for
\eqref{RightWhittEqs} and $\frak{n}^{(m,N)}_+$. $\Box$

\begin{de}
Given a pair of characters $\chi_{\pm}$ of the opposed nilpotent
subalgebras $\frak{n}^{(m,N)}_{\pm}$, the $\Gr_{m,N}$-Whittaker
function associated with the principal series representation
$\bigl(\pi_{\ul},\,\CV_{\ul}\bigr)$ is defined as the following
matrix element:
 \be\label{GrWhittaker}
  \Psi^{(m,N)}_{\underline{\lambda}}(\ux)\,
  =\,e^{-x_1\frac{m(N-m)}{2}}
  \bigl\<\psi_L^{(m,N)}\,,\,
  \pi_{\underline{\lambda}}\bigl(g(x_1,\ldots,x_N)\bigr)\,
  \psi_R^{(m,N)}\bigr\>\,,
 \ee
where  the left and right  vectors  solve the equations
\eqref{LeftWhittEqs} and \eqref{RightWhittEqs} respectively, and .
Here $g(x)$ is a $H^{(m)}$-valued function given by
 \be\label{GrGroupElement}
  g(\ux)\,=\,\exp\Big\{\sum_{i=1}^N\,x_iH_i\Big\}\,,
 \ee
where $\ux=(x_1,\ldots,x_N)$ and the generators $H_i$,
$i=1,\ldots,N$ are defined by \eqref{ModTriDecompositionH}.
\end{de}
In the above definition $\<\,,\,\>$ denotes a non-degenerate pairing
between the Whittaker submodules: $\CW_L\times\CW_R\to\IC$.

\subsection{Stationary phase integral
and Toric degeneration of $\Gr_{m,N}$}

This part contains our main result; namely, we introduce the
stationary phase integral representation for the
$\Gr_{m,N}$-Whittaker function, and then establish a direct
connection to the toric degeneration of Grassmannian $\Gr_{m,N}$
proposed in \cite{BCFKvS1}.
\begin{te}
The specialized $\Gr_{m,N}$-Whittaker function \eqref{GrWhittaker}
has the following stationary phase integral representation.
 \be\label{GivIntegralRep}
  \Psi^{(m,N)}_{\underline{\lambda}}(x_{N,1},0,\ldots,0)\,
  =\,\int\limits_{\CC}\,\omega_{m,\,N}\quad
  e^{^{\CF_{m,N}(\ux)}}\,,
 \ee
where
 \be\label{GrGivAction}
  \CF_{m,N}(\ux)\,\,
  =\,\,\i\Big(\sum_{k=1}^m\la_{N-m+k}\Big)x_{N,1}\,
  +\,\i\sum_{n=1}^{N-m}(\la_n-\la_{n+1})\!\!
  \sum_{i=1}^{\min(m,n)}\!\!x_{n,i}\\
  +\,\,\i\sum_{n=1}^{m-1}(\la_{N-m+n}-\la_{N-m+n+1})
  \sum_{i=n+1}^{\min(N-m+n,\,m)}x_{N-m+n,\,i}\\
  -\,\frac{1}{\h}\Big(e^{-x_{mm}}\,
  +\,e^{x_{N-m,\,1}-x_{N,1}}\,
  +\,\sum_{k=1}^m\sum_{i=1}^{N-1-m}e^{x_{i+k-1,\,k}-x_{i+k,\,k}}\\
  +\,\sum_{k=1}^{N-m}\sum_{i=1}^{m-1}e^{x_{k+i,\,i+1}-x_{k+i-1,\,i}}
  \Big)\,.
 \ee
and
 \be\label{FlatMeasure}
  \omega_{m,N}\,
  =\,\prod_{n=1}^{N-m}\prod_{k=1}^{\min(n,\,m)}dx_{n,k}\cdot
  \prod_{n=1}^{m-1}\,\,\prod_{i=n+1}^{\min(N-m+n,\,m)}dx_{N-m+n,\,i}\,.
 \ee
The integration contour $\CC$ is a slight deformation of
$\IR^{m(N-m)}$ in $\IC^{m(N-m)}$ such that the integrand decreases
exponentially.
\end{te}
\emph{Proof:} The proof is given in Section 3. $\Box$

Specifying the parameters $\la_n=0,\,n=1,\ldots,N$, the function
$\CF_{m,N}(\ux)$ admits a simple combinatorial structure. Namely,
let us consider the following graph:
 \be\label{GrGraph}
 \xymatrix{
 x_{N,\,1}\ar@/_1pc/[dd] &&&\\
 \vdots & \ddots & &\\
 x_{N-m,\,1}\ar[r]\ar[d] & \ldots\ar[r]\ar[d] & x_{N-1,\,m}\ar[d] & \\
 \vdots\ar[d]\ar[r] & \ddots\ar[r]\ar[d] & \vdots\ar[d] & \\
 x_{21}\ar[d]\ar[r] & \ddots\ar[r]\ar[d] & x_{m+1,\,m}\ar[d] & \\
 x_{11}\ar[r] & \ldots\ar[r] & x_{m,\,m}\ar[r] & 0
 }
 \ee
Then let us associate with every arrow $x\to y$ an exponential
function $e^{y-x}$; to any interior vertex $x_{k,i}$ in
\eqref{GrGraph} let us assign a pair of exponential functions
$a_{k,i}=e^{x_{k-1}-x_{k,i}}$ and
$b_{k,i}=e^{x_{k+1,\,i+1}-x_{k,i}}$. Besides, let
$a_N=e^{x_{N-m,\,1}-x_{N,1}}$ and $b_m=e^{-x_{mm}}$.
\begin{cor} The function $\CF_{m,N}(\ux)$ \eqref{GrGivAction} equals to
the sum of exponential functions for all the arrows in the graph
\eqref{GrGraph}:
$$
 \CF_{m,N}\,=\,a_N+\sum_{k=2}^{N-m}a_{k,\,m}\,
 +\,b_m+\sum_{i=1}^{m-1}b_{1,\,i}\,
 +\,\sum_{k=2}^{N-m}\sum_{i=1}^{m-1}a_{k,i}+b_{k,i}\,.
$$
\end{cor}

Actually, the graph \eqref{GrGraph} defines a toric degeneration of
the Grassmann variety \cite{BCFKvS1}, \cite{BCFKvS2}. Namely, the
torification of  $\Gr_{m,N}$ can be identified with a spectrum of
the algebra of functions in $a_{ki},\,k=2,\ldots,N-m,\,i=1,\ldots,m$
and $b_{nj},\,n=1,\ldots,N-m,\,j=1,\ldots,m-1$ modulo the ideal of
relations:
 \be\label{ToricVariety}
  a_{k,\,i}b_{k-1,\,i}\,=\,b_{k,\,i}a_{k+1,\,i+1}\,,
 \hspace{1.5cm}
  k=2,\ldots,N-m,\quad i=1,\ldots,m-1\,;\\
  a_Nb_m\prod_{i=1}^{m-1}a_{N-m,\,i}\prod_{k=2}^{N-m}b_{k,\,m}\,\,
  =\,\,q\,,
 \hspace{2.5cm}
  q=e^{-x_1}\,.
 \ee

\subsection{Quantum $\Gr_{m,\,N}$-Toda lattice}

Actually, the $\Gr_{m,N}$-Whittaker function defines a D-module
introduced in \cite{GLO2} and called quantum $\Gr_{m,\,N}$-Toda
lattice. This D-module is provided by the infinitesimal action of
the universal enveloping algebra $\CU(\gl_N)$ in our representation
$(\pi_{\ul},\,\CV_{\ul})$. In this part we describe the D-module
$\CD_{m,N}$ defined by \eqref{GrWhittaker}, and then identify
$\CD_{m,N}$ with the $\Gr_{m,\,N}$-Toda lattice from \cite{GLO2}.

The action of the center $\CZ\subset\CU(\gl_N)$ of the universal
enveloping algebra in principal series representation
$(\pi_{\ul},\,\CV_{\ul})$ produces the following action of
differential operators, the parabolic Toda lattice Hamiltonians, on
$\Gr_{m,N}$-Whittaker function:
 \be\label{GrTodaHamiltonians}
  \CH^{(m,N)}_k(x,\pr_x)\cdot\Psi_{\ul}^{(m,\,N)}(\ux)\,
  =\,\hbar^k\,e^{-x_1\frac{m(N-m)}{2}}\bigl
  \<\psi_L\,,\,\pi_{\underline{\la}}\bigl(c_k\,g(\ux)\bigr)
  \psi_R\bigr\>\,,
 \ee
for $k=1,\ldots,N$ with $c_k,\,k=1,\ldots,N$ being the Casimir
generators of the center $\CZ$. The first two Casimir elements are
given by:
 \be\label{glCasimirs}
  C_1\,=\,\sum_{i=1}^NE_{ii},\,
 \hspace{1.5cm}
  C_2\,=\,\sum_{i,j=1\atop i<j}^N\Big(E_{ii}E_{jj}
  -E_{ji}E_{ij}+\rho_i\rho_j\,\Big)\,-\,\sum_{i=1}^N\rho_iE_{ii}\,,
 \ee
where $\rho_i=(N+1-2i)/2,\,i=1,\ldots,N$.
\begin{prop} Action of the first two Casimir generators
\eqref{glCasimirs} have the following explicit form:
 \be\label{GrHamiltonians}
  \CH_1^{(m,N)}\,
  =\,\h\frac{\pr}{\pr x_1}\,+\,\h\frac{\pr}{\pr x_N}\,,\\
  \CH_2^{(m,N)}\,
  =\,\h^2\Big\{\frac{\pr^2}{\pr x_1\pr x_N}\,
  +\,\sum_{1\leq i\leq j\leq m}\Big(\,(-x_i)^{1-\delta_{i,\,1}}
  \frac{\pr}{\pr x_i}\Big)\Big(x_j\frac{\pr}{\pr x_j}\Big)\\
  +\,\sum_{m+1\leq i\leq j\leq N}\Big(x_i\frac{\pr}{\pr x_i}\Big)
  \Big(\,x_j^{1-\delta_{j,\,N}}\frac{\pr}{\pr x_j}\Big)\,
  -\,\sum_{k=1}^m(k-1)x_k\frac{\pr}{\pr x_k}\\
  -\!\sum_{k=m+1}^N(N+1-k)x_k\frac{\pr}{\pr x_k}\Big\}\,
  -\,\h\Big\{\sum_{i=1}^{m-1}(-x_i)^{1-\delta_{i,\,1}}
  \frac{\pr}{\pr x_{i+1}}\,
  +\!\sum_{j=m+1}^{N-1}x_{j+1}^{1-\delta_{j,\,N-1}}
  \frac{\pr}{\pr x_j}\Big\}\\
  +\,(-1)^{\delta_{m,\,N-1}+\e(m,N)}x_m^{1-\delta_{m,1}}
  x_{m+1}^{1-\delta_{m,\,N-1}}\,e^{x_N-x_1}\,
  -\,\frac{\h^2}{24}(N-1)(N-2)(N-3)\,.
 \ee
\end{prop}
\emph{Proof:} The first statement is trivial. The proof of the
second formula is given in Section 4. $\Box$
\begin{rem}
The Hamiltonians \eqref{GrHamiltonians} coincide (up to signs in
accordance with the choice of signs in \eqref{GrGroupElement}) with
the first two $\Gr_{m,N}$-Toda Hamiltonians from \cite{GLO2}. Let us
emphasize that although for generic $m$ the $\Gr_{m,N}$-Whittaker
vectors \eqref{LeftWhittEqs} and \eqref{RightWhittEqs} are different
from the ones introduced in \cite{GLO2}, the two (generalized)
Whittaker models: the one from \cite{GLO2}, and its modification
introduced above, produce the same $\Gr_{m,N}$-Toda D-module. In
particular, the Hamiltonians \eqref{GrTodaHamiltonians} are
identical to the parabolic Toda Hamiltonians from \cite{GLO2}, and
after specifying $x_2=\ldots=x_N=0$ the symbols of Hamiltonians
$\CH_k^{(m,N)},\,k=1,\ldots,N$ generate the small quantum cohomology
algebra $qH^*(\Gr_{m,N})$.
\end{rem}

Let $\CD_{m,N}$ be the D-module generated by the Hamiltonians
$\CH_k^{(m,N)},\,k=1,\ldots,N$ as a module over the algebra of
differential operators $\CD\bigl(e^x,\pr_x,\h\bigr)$:
$$
 \CD_{m,N}\,\simeq\,\CD\bigl(e^x,\pr_x,\h\bigr)\bigr/\CT_{m,N}\,,
 \hspace{1.5cm}
 \CT_{m,N}\,=\,\bigl\<
 \CH^{(m,N)}_1(x,\pr_x),\,\ldots,\,\CH^{(m,N)}_N(x,\pr_x)\,\bigr\>\,.
$$
Equivalently, the D-module $\CD_{m,N}$ can be defined in terms of
the quantum Lax operator, the following $\Mat(N,\IR)$-valued
differential operator:
 \be\label{Loper}
  \CL(\ux,\pr_{\ux})\,\cdot\Psi_{\ul}^{(m,\,N)}(\ux)\,:=\,
  \h\sum_{i,j=1}^N\,e_{ij}\,e^{-x_1\frac{m(N-m)}{2}}\,
  \<\psi_L\,,\,\pi_{\ul}\bigl(E_{ij}\, g(\ux)\bigr)\psi_R\>,
 \ee
where $(e_{ij})_{kn}=\delta_{ik}\delta_{jn},\,\,$ for
$i,j,k,n=1,\ldots ,N$ are the matrix units.
\begin{prop}
The quantum Lax operator
$\CL(x_1,\ldots,x_N;\,\pr_{x_1},\ldots,\pr_{x_N})
=\|\CL_{ij}\|,\,i,j=1,\ldots,N$
has the following form:
  \be\label{GrLaxOperator}
   \CL_{k,1}\,=\,\h\pr_{x_k}\,,\quad k=1,\ldots,m\,;
  \hspace{1.5cm}
   \CL_{m+1,\,1}\,=\,-1\,,\\
   \CL_{k,1}\,=\,0\,,\quad k=m+2,\ldots,N\,;
  \hspace{1.2cm}
   \CL_{k,j}\,=\,0\,,\quad j=2,\ldots,m,\,\, k=j,\ldots,N\,;\\
   \CL_{a+1,\,a}\,=\,-1\,,
  \hspace{1.2cm}
   \CL_{k,\,a}\,=\,0\,,\quad k=a+2,\ldots,N,\,\,a=m+1,\ldots,N-1\,;\\
   \CL_{1,\,k}\,=\,-\delta_{k,2}+(1-\delta_{k,\,m})x_{k+1}\,
   +\,x_k\pr_{x_1}\,+\,\sum_{n=2}^mx_kx_n\pr_{x_n}\,,\qquad
   k=2,\ldots,m\,;\\
   \CL_{k,\,i}\,=\,\delta_{i,\,k+1}\,+\,\h x_i\pr_{x_k}\,,\qquad
   k=2,\ldots,m-1\,,\quad i=k+1,\ldots,m\,;\\
   \CL_{1,\,a}\,=\,-(-1)^{\e(m,N)}x_ax_me^{x_N-x_1}\,,
  \hspace{1cm}
   \CL_{k,\,a}\,=\,0\,,\quad k=2,\ldots,m-1\,,\\
   \CL_{m,\,a}\,=\,(-1)^{\e(m,N)}x_ae^{x_N-x_1}\,,
   \qquad a=m+1,\ldots,N-1\,;\\
   \CL_{m,\,N}\,=\,-(-1)^{\e(m,N)}e^{x_N-x_1}\,;\\
   \CL_{aa}\,=\,\h x_a\pr_{x_a}\,,
  \hspace{1.2cm}
   \CL_{a,\,N}\,=\,\h\pr_{x_a}\,,\quad a=m+1,\ldots,N-1\,;\\
   \CL_{NN}\,=\,\h\pr_{x_N}\,-\,\sum_{a=m+1}^{N-1}x_a\pr_{x_a}\,.
  \ee
\end{prop}
\emph{Proof:} The proof is given in Section 4. $\Box$

\noindent The symbol of the quantum Lax operator is referred to as
the Lax matrix $L(x_1,\ldots,x_N;\,p_1,\ldots,p_N)$.
\begin{cor}
(i) The Lax matrix $L(\ux;\,\upp)$ of \eqref{GrLaxOperator} and the
Lax matrix introduced in \cite{GLO2} have identical characteristic
polynomials.

(ii) The specialized Lax matrix
$L(x_1,0,\ldots,0;\,p_1,\ldots,p_N)=\|L_{ij}\|,\,i,j=1,\ldots,N$ is
given by
 \be\label{GrLaxMatrix}
  L_{k,1}\,=\,p_k\,,\quad k=1,\ldots,m\,;
 \hspace{1.5cm}
   L_{m+1,\,1}\,=\,-1\,,\\
   L_{k,1}\,=\,0\,,\quad k=m+2,\ldots,N\,;
  \hspace{1.2cm}
   L_{k,j}\,=\,0\,,\quad j=2,\ldots,m,\,\, k=j,\ldots,N\,;\\
   L_{a+1,\,a}\,=\,-1\,,\quad
   L_{aa}\,=\,L_{k,a}\,=\,0\,,\quad k=a+2,\ldots,N,\quad
   a=m+1,\ldots,N-1\,;\\
   L_{i,\,i+1}\,=\,-1\,,\quad i=1,\ldots,m-1\,;\\
   L_{a,N}\,=\,-p_a\,,\quad a=m+1,\ldots,N-1\,;
  \hspace{1.5cm}
   L_{m,N}\,=\,-(-1)^{\e(m,N)}e^{-x_1}
 \ee
It defines the small quantum cohomology algebra $qH^*(\Gr_{m,N})$ in
the form of \cite{AS}, \cite{K}.
\end{cor}
\emph{Proof:} One can readily check that the matrix $\|L_{ij}\|$ and
matrix $A$ in \cite{AS} defining $qH^*(\Gr_{m,N})$ have identical
characteristic polynomials, and thus
$\det\bigl(\la+\|L_{ij}\|\bigr)$ is the generating function of the
ideal for the small quantum cohomology algebra. $\Box$

\section{Gelfand-Zetlin graphs, paths, and $\Gr_{m,N}$-Whittaker vectors}

In this section we recall the Gauss-Givental realization of the
universal enveloping algebra $\CU=\CU(\gl_N)$ introduced in
\cite{GKLO}. In the second part of this Section we apply this
construction to derivation of the $\Gr_{m,N}$-Whittaker vectors,
solving the defining relations \eqref{LeftWhittEqs} and
\eqref{RightWhittEqs}.

Actually, the construction of Gauss-Givental realization of
principal series $\CU$-modules originates from the total positivity
phenomenon in unipotent varieties developed by Lusztig \cite{L}; a
detailed study of connections between Gauss-Givental realizations of
$\CU(\frak{g})$ and total positivity can be found in \cite{GLO1}.

\npaa Let $\CM_N$ be the space of meromorphic functions in
$e^{x_{n,k}},\,n=1,\ldots,N-1\,;\,k=1,\ldots,n$, then the standard
generators $E_{ij},\,i,j=1,\ldots,N$ of $\gl_N$ admit the following
realization by first-order differential operators in $\CM_N$:
 \bqa
  E_{i,i}&=& \mu_i\,-\,
  \sum_{k=1}^{i-1}\frac{\pr}{\pr x_{N+k-i,k}}\,
  +\,\sum_{k=i}^{N-1}\frac{\pr}{\pr x_{k,i}},\nonumber \\
  E_{i,i+1}&=& -\sum_{n=1}^ie^{x_{N-1-i+n,\,n}-x_{N-i+n,\,n}}
  \sum_{k=1}^n\Big\{\frac{\pr}{\pr x_{N-1-i+k,\,k}}\,
  -\,\frac{\pr}{\pr x_{N-1-i+k,\,k-1}}\Big\}\,,
 \label{GGivRep}\\
  E_{i+1,i}&=&\sum_{n=1}^{N-i}e^{x_{n+i,\,i+1}-x_{k+i-1,\,i}}
  \Big[\mu_i-\mu_{i+1}\,+\,\sum_{k=1}^n\Big\{
  \frac{\pr}{\pr x_{i+k-1,\,i}}\,-\,\frac{\pr}{\pr x_{i+k-1,\,i+1}}
  \Big\}\Big]\,,
  \nonumber
 \eqa
where $x_{N,\,i}=0,\, i=1,\ldots,N$ is assumed.

The universal enveloping algebra acts in $\CV_{\ul}\subset\CM_N$ by
differential operators \eqref{GGivRep} with
 \be\label{SpectralParameters}
  \mu_n\,=\,\i\la_n\,-\,\rho_n^{(N)}\,,
 \hspace{1cm}
  \rho_n\,=\,n-\frac{N+1}{2}\,,
 \hspace{1.5cm}
  n=1,\ldots,N\,,
 \ee
and the Whittaker submodules
$\CW_{L,\,R}\subset\CV_{\ul}\subset\CM_N$ are spanned by
$\prod\limits_{1\leq k\leq i\leq N-1}e^{n_{k,i}x_{k,i}}\psi_{L,\,R}$
with $n_{k,i}\in\IZ$. The non-degenerate pairing between the
Whittaker modules is given by
 \be\label{GGivPairing}
  \<\phi_1,\,\phi_2\>\,=\,\int_{\CC}\mu_N(x)\,\ov{\phi_1}\,\phi_2\,,
 \hspace{1.5cm}
  \phi_1\in\CW_L,\,\quad \phi_2\in\CW_R\,,
 \ee
where the integration contour is a slight deformation of
$\IR^{N(N-1)/2}$ in $\IC^{N(N-1)/2}$ such that the integrand
exponentially decreases for $\phi_1=\psi_L$ and $\phi_2=\psi_R$, and
the measure $\mu_N(x)$ is given by
 \be\label{GGivMeasure}
  \mu_N(x)\,=\,\prod_{k=1}^{N-1}\prod_{i=1}^ke^{-x_{k,\,i}}dx_{k,\,i}\,.
 \ee
One can readily check that thus defined pairing \eqref{GGivPairing}
between $\CW_L$ and $\CW_R$ possesses the following property:
 \be\label{HermiteanProp}
  \<X\cdot\phi_1,\,\phi_2\>\,=\,-\<\phi_1,\,X\cdot\psi_2\>\,,
 \hspace{1.5cm}
  X\in\gl_N\,,\quad
  \phi_1\in\CW_L,\,\quad \phi_2\in\CW_R\,.
 \ee

\npaa The Gauss-Givental realization of $\CU(\gl_N)$ possesses a
distinguished combinatorial structure arising from the
Gelfand-Zetlin graph (see \cite{GLO1}):
 \bqa\label{GZgraph}
  \xymatrix{
  x_{N,\,1}\ar[d] &&&&\\
  x_{N-1,1}\ar[d]\ar[r] & x_{N,\,2}\ar[d] &&&\\
  \vdots\ar[d] & \ddots & \ddots &&\\
  x_{21}\ar[d]\ar[r] & \ddots\ar[d] & \ddots\ar[r] &
  x_{N,\,N-1}\ar[d]\\
  x_{11}\ar[r] & x_{22}\ar[r] &
  \ldots\ar[r] & x_{N-1,\,N-1}\ar[r] & x_{NN}
  }
 \eqa
Namely, let $\CI_N$ be the set of vertices in the Gelfand-Zetlin
graph:
$$
 \CI_N\,
 =\,\bigl\{(n,\,j)\in\IZ_+^2\,;\,1\leq j\leq n\leq N\bigr\}\,;
$$
there is a tautological embedding $\CI_m\subset\CI_N$, for any
$0<m<N$. Given $(n,j)\in\CI_{N-1}\subset\CI_N$ let $A^r_{n,\,j}$ be
the following function attached to a vertex $x_{n,\,i}$ in
\eqref{GZgraph}:
 \be\label{PathFunctionA}
  A_{n,\,j}^r\,\,
  =\,\,\sum_{I_r}\,\prod_{\a=1}^r\,
  e^{x_{n+i_{\a}-\a,\,j+i_{\a}}-x_{n+i_{\a}-\a-1,\,j+i_{\a}}}\,,
 \ee
where the summation goes over the strict partitions
$$
 I_r\,=\,(i_1<\ldots<i_r)\in\IZ_+^r\,,
\hspace{1.5cm}
 i_{\a}\leq N-n+\a,\quad\a=1,\ldots,r\,.
$$
The function $A_{n,\,j}^r(k)$ satisfy the following evident
relation:
 \be\label{PathRelationA}
  A^r_{n,j}\,
  =\,A^r_{n+1,\,j+1}\,
  +\,A^{r-1}_{n,\,j+1}\,e^{x_{n,\,j+1}-x_{n+1,\,j+1}}\,.
 \ee

Also for $(n,\,j)\in\CI_{N-1}\subset\CI_N$ let us introduce the
function $B_{n,\,i}$ given by
 \be\label{PathFunctionB}
  B_{n,\,j}(k)\,\,
  =\,\,\sum_{I^*_{k+j-n-1}}\,e^{x_{n,\,j}-x_{n+1,\,j}}\,
  \prod_{\a=1}^{k+j-n-1}\,
  e^{x_{n+i_{\a}+\a,\,j+\i_{\a}}-x_{n+i_{\a}-\a+1,\,j+i_{\a}}}\,,
 \ee
where the summation goes over the partitions
$$
 I^*_{k+j-n-1}\,=\,(i_1\leq\ldots\leq i_{k+j-n-1})\in\IZ_+^{k+j-n-1}\,.
$$
The function $B_{n,\,j}$ satisfies the following relation:
 \be\label{PathRelationB}
  B_{n,\,j}(k)\,e^{x_{n+1,\,j}-x_{n,\,j}}\,\,
  =\,\,B_{n+1,\,j}(k)\,\,+\,\,B_{n+2,\,j+1}(k)\,.
 \ee
Actually, the relations \eqref{PathRelationA} and
\eqref{PathRelationB} are direct consequence of the "box relations"
\eqref{ToricVariety}.

\npaa With respect to an obvious symmetry of the graph
\eqref{GrGraph}, let us introduce the following pair of functions:
 \be\label{PathFunctionAt}
  \At^r_{n,\,j}\,\,=\,\,\sum_{I^*_{r-1}}\,e^{x_{n+1,\,j+1}-x_{n,\,j}}\,
  \prod_{\a=1}^{r-1}\,
  e^{x_{n+i_{\a}+\a+1,\,j+\a+1}-x_{n+i_{\a}\a,\,j+\a}}\,,
 \ee
where the summation goes over partitions
$$
 I^*_{r-1}\,=\,(i_1\leq\ldots\leq i_{r-1})\in\IZ_{\geq0}^{r-1}\,,
\hspace{1.5cm}
 i_{\a}\leq N-n-r\,,\quad \a=1,\ldots,r-1\,;
$$
and
 \be\label{PathFunctionBt}
  \Bt_{n,\,j}(k)\,\,=\,\,\sum_{I_{j-k}}\,\,\prod_{\a=1}^{j-k}\,
  e^{x_{n+i_{\a}-\a+1,\,j+1-\a}-x_{n+i_{\a}-\a,\,j-\a}}\,,
 \ee
where the summation goes over strict partitions
$$
 I_{j-k}\,=\,(i_1<\ldots<i_{j-k})\in\IZ_+^{j-k}\,,
\hspace{1.5cm}
 i_{\a}\leq N-n+\a\,,\quad \a=1,\ldots,j-1\,.
$$
Analogously to the the functions $A^r_{n,\,j}$ and $B_{n,\,j}$, the
"box relations" \eqref{ToricVariety} imply the following relations
for the introduced functions \eqref{PathFunctionAt} and
\eqref{PathFunctionBt}:
 \be\label{PathRelationAt}
  \At^r_{n,\,j}\,e^{x_{n+1,\,j+1}-x_{n,\,j}}\,\,
  =\,\,\At^{r-1}_{n+1,\,j+1}\,+\,\At^{r-1}_{n+2,\,j+1}\,;
 \ee
and
 \be\label{PathRelationBt}
  \Bt_{n,\,j}\,\,
  =\,\,\Bt_{n+1,\,j}\,+\,\Bt_{n,\,j-1}\,e^{x_{n+1,\,j}-x_{n,\,j-1}}\,.
 \ee

\npaa In fact, the summations in \eqref{PathFunctionA},
\eqref{PathRelationAt}, and \eqref{PathFunctionB},
\eqref{PathRelationBt} can be readily interpreted as sums over paths
in Gelfand-Zetlin graph. More precisely, the functions
 \be\label{PathFunctions}
  P_{n,\,j}^r(k)\,\,:=\,\,A^{r+n-k-j}_{n,\,j}\cdot B_{n,\,j}(k)\,,
 \hspace{1.5cm}
  \Pt_{n,\,j}^r\,\,:=\,\,\At^{r+k-j}_{n,\,j}\cdot\Bt_{n,\,j}\,,
 \ee
are represent sums over all paths (with certain restrictions) of
length $r$ passing through a vertex $x_{n,\,j}$ on the
Gelfand-Zetlin graph; the paths from $P_{n,\,j}^r(k)$ are starting
at horizontal line $\{x_{a,\,b};\,a-b=k\}$, and the paths from
$\Pt^r_{n,\,j}(k)$ are starting at vertical line
$\{x_{a,\,b};\,b=k\}$ in graph \eqref{GZgraph}.

Moreover, the generators $E_{ij}$ of Lie algebra $\gl_N$ in
Gauss-Givental realization \eqref{GGivRep} admit a distinct
description in terms of certain paths in graph \eqref{GZgraph}.
Namely, for any vertex $x_{n,\,j}$ let us introduce the following
pair of differential operators:
 \be
  D_{n,\,j}\,\,=\,\,\frac{\pr}{\pr x_{n+1-j,\,1}}\,
  +\,\sum_{i=1}^{j-1}\Big(\frac{\pr}{\pr x_{n+1+i-j,\,i+1}}
  -\frac{\pr}{\pr x_{n+1+i-j,\,i}}\Big)\,,\\
  \Dt^{\mu}_{n,\,j}\,\,=\,\,\mu_j-\mu_{j+1}\,
  +\,\frac{\pr}{\pr x_{j,\,j}}\,
  +\,\sum_{i=1}^{n-j}\Big(\frac{\pr}{\pr x_{i+j,\,j}}
  -\frac{\pr}{\pr x_{i+j,\,j+1}}\Big)\,.
 \ee
\begin{prop} The Lie algebra generators $E_{ij}$ have the following
combinatorial realization in terms of the Gelfand-Zetlin graph
\eqref{GZgraph}:
 \be
  E_{n,\,j}\,=\,\sum_{k=j}^{n-1}(-1)^{k+1}\sum_{i=0}^{N-n}
  \Pt^{n-j}_{k+i,\,k}(j)\,\Dt^{\mu}_{k+i,\,k}\,,
 \hspace{1.5cm}
  n>j\,;\\
  E_{n,\,i}\,=\,\sum_{k=1}^{i-n}(-1)^k\sum_{j=0}^{n-1}
  P^{i-n}_{N+k-i+j,\,j+1}(N-n)\,D_{N-i+k+j,\,j+1}\,,
 \ee
for $n<i$.
\end{prop}
\begin{cor} The elements $E_{m,\,1}$ and $E_{m,N}$ have the following
differential operators in the Gauss-Givental realization
\eqref{GGivRep}.
 \be\label{EkN}
 \hspace{-1.5cm}
  E_{n,\,N}\,=\,\sum_{k=1}^{N+1-n}(-1)^{k+1}\sum_{i=1}^n\,
  P^{n+1-i}_{k+i-1,\,i}\,\Big\{
  \frac{\pr}{\pr x_{k+i-1,\,1}}\,
  +\,\sum_{j=2}^i\Big(\frac{\pr}{\pr x_{k+i-2,\,i}}\,
  -\,\frac{\pr}{\pr x_{k+i-2,\,i}}\Big)\Big\}\,;\\
 \hspace{-1.5cm}
  E_{n,\,1}\,=\,\sum_{k=1}^{n-1}(-1)^k\sum_{i=1}^{N-n}\,
  \Pt^{n-1}_{k+i,\,i}\Big\{\mu_k-\mu_{k+1}\,
  +\,\frac{\pr}{\pr x_{kk}}\,
  +\,\sum_{j=1}^i\Big(\frac{\pr}{\pr x_{k+j,\,k}}
  -\frac{\pr}{\pr x_{k+j,\,k+1}}\Big)\Big\}\,,
 \ee
where $P^r_{n,\,j}:=P^r_{n,\,j}(1)$, and
$\Pt^r_{}:=\Pt^r_{n,\,j}(N-1)$.
\end{cor}
\emph{Proof:} Direct calculations of commutators, using
\eqref{GGivRep}. $\Box$

\subsection{Derivation of $\Gr_{m,N}$-Whittaker vectors}

In this part we solve the defining equations \eqref{LeftWhittEqs}
and \eqref{RightWhittEqs} and find out the $\Gr_{m,N}$-Whittaker
vectors.
\begin{prop} The following $\Gr_{m,N}$-Whittaker vectors satisfy the defining
equations \eqref{LeftWhittEqs} and \eqref{RightWhittEqs} in
realization \eqref{GGivRep}:
 \be\label{GrLeftWhittaker}
  \psi_L^{(m,N)}\,=\,\frac{1}{C^L_{m,N}}\,
  \exp\Big\{-\sum_{n=1}^{N-1}\bigl(\mu_n-\mu_{n+1}\bigr)\sum_{i=1}^nx_{n,i}\,
  +\,\sum_{k=1}^{m-1}\mu_{N-1-k}x_{N-1-k,\,1}\,\\
  -\,\frac{1}{\h}\Big(e^{x_{N-m,\,1}}\,
  +\,\sum_{k=1}^m\sum_{i=1}^{N-m-1}e^{x_{i+k-1,\,k}-x_{i+k,\,k}}\,
  +\,\sum_{k=m+1}^{N-1}\Big[e^{x_{N-1,\,k}}\\
  +\,\sum_{i=1}^{N-k-1}e^{x_{i+k-1,\,k}-x_{i+k,\,k}}\Big]\Big)\Big\}
 \ee
and
 \be\label{GrRightWhittaker}
  \psi_R^{(m,N)}\,=\,\frac{1}{C^R_{m,N}}\,
  \exp\Big\{-\sum_{k=m+1}^{N-1}\mu_kx_{kk}\,
  -\,\frac{1}{\h}\Big(e^{-x_{mm}}\,
  +\,\sum_{k=1}^{N-m}\sum_{i=1}^{m-1}e^{x_{k+i,\,i+1}-x_{k+i-1,\,i}}\\
  +\,\sum_{k=1}^{m-1}\Big[e^{-x_{N-1,\,k}}\,
  +\,\sum_{i=1}^{k-1}e^{x_{N-k+i,\,i+1}-x_{N-k+i-1,\,i}}\Big]\Big)\Big\}
 \ee
where
 \be\label{NormalizationCfunctions}
  C^L_{m,N}\,=\,\prod_{i,j=N-m+1\atop i<j}^N\h^{\rho_i-\i\la_j}
  \Gamma\bigl(-\i\la_j-\rho_i\bigr)\,,
 \hspace{1cm}
  C^R_{m,N}\,=\,\prod_{i,j=m+1\atop i<j}^N\h^{\i\la_j-\rho_i}
  \Gamma\bigl(\rho_i-\i\la_j\bigr)\,.
 \ee
\end{prop}
\emph{Proof.} Our proof of Proposition 2.1 is based on an
verification of the defining equations \eqref{LeftWhittEqs},
\eqref{RightWhittEqs}, using the Gauss-Givental realization
\eqref{GGivRep}.

Actually, the expressions \eqref{GrLeftWhittaker} and
\eqref{GrRightWhittaker} (with specialized parameters
$\mu_n=0,\,n=1,\ldots,N$) have definite interpretation in terms of
arrows in \eqref{GZgraph} defined by the equations
\eqref{LeftWhittEqs}, \eqref{RightWhittEqs}, respectively. In this
way the graph \eqref{GrGraph} is a subgraph of \eqref{GZgraph},
built of the corresponding arrows from \eqref{GrLeftWhittaker} and
\eqref{GrRightWhittaker}.

\npaa At first let us observe that the action of the Cartan
generators $E_{ii},\,i=2,\ldots,N-1$ fixes a dependence of the
Whittaker vectors on the parameters $\mu_1,\ldots,\mu_N)$. Namely,
the following holds:
$$
 E_{kk}\cdot
 \exp\Big\{-\sum_{n=1}^{N-1}\bigl(\mu_n-\mu_{n+1}\bigr)\sum_{i=1}^nx_{n,i}\,
 +\,\sum_{k=1}^{m-1}\mu_{N-1-k}x_{N-1-k,\,1}\Big\}\,=\,0
$$
for $k=2,\ldots,m$, and
$$
 E_{aa}\cdot e^{-(\mu_{m+1}x_{m+1,\,m+1}+\ldots+\mu_{N-1}x_{N-1,\,N-1})}\,=\,0\,,
\hspace{1.5cm}
 a=m+1,\ldots,N-1\,.
$$

\npaa Besides, the action of differential operators $E_{ii}$ have
the following properties:
$$
 E_{kk}\cdot e^{x_{n,k}-x_{n+1,\,k}}\,=\,0\,,
\hspace{1.5cm}
 k=2,\ldots,m\,,
$$
when $n=2,\ldots,N-m,\,k=2,\ldots,m$, and
$$
 E_{aa}\cdot e^{x_{n+1,\,k+1}-x_{n,k}}\,=\,0\,,
\hspace{1.5cm}
 a=m+1,\dots,N-1\,,
$$
when $n=2,\ldots,N+1-m,\,k=1,\ldots,m-1$. Also, taking into
account that $E_{kk},\,k=2,\ldots,m$ annihilate any function in
$x_{n,1},\,n=1,\ldots,N-m$, and $E_{aa},\,m+1,\ldots,N-1$ annihilate
any function in $x_{nn},\,n=1,\ldots,m$, one can deduce that
\eqref{GrLeftWhittaker} and \eqref{GrRightWhittaker} satisfy
$$
 E_{kk}\psi_L^{(m,N)}\,=\,0\,,\quad k=2,\ldots,m\,;
\hspace{1.5cm}
 E_{aa}\psi_R^{(m,N)}\,=\,0\,,\quad a=m+1,\ldots,N-1\,.
$$

\npaa Next for $\psi_R^{(m,N)}$ one finds out the following:
 \be
  E_{k-1,\,k}\psi_R^{(m,N)}\,=\,-\frac{1}{\h}\psi_R^{(m,N)}\,,
 \hspace{1.5cm}
  k=2,\ldots,m\,,
 \ee
since
$$
 E_{i,\,i+1}\cdot e^{-\h^{-1}\,e^{-x_{N-1,\,i}}}\,
 =\,-\frac{1}{\h}\,e^{-\h^{-1}\,e^{-x_{N-1,\,i}}}\,,
\hspace{1.5cm}
 i=1,\ldots,N-1\,,
$$
and
$$
 E_{i,\,i+1}\cdot\exp\Big\{-\frac{1}{\h}\sum_{j=1}^{m-1}
 e^{x_{N+j-i,\,j+1}-x_{N+j-i-1,\,j}}\,
 +\,e^{x_{N+j-i+1,\,j+1}-x_{N+j-i,\,j}}
 \Big\}\,=\,0\,,\quad
 i=1,\ldots,N-1
$$
for any $m=1,\ldots,N-1$, is due to the "box relations"
(\ref{ToricVariety}.1):
$$
 e^{x_{N+j-i-1,\,j}-x_{N+j-i,\,j}}\cdot
 e^{x_{N+j-i,\,j+1}-x_{N+j-i-1,\,j}}\,
 =\,e^{x_{N+j-i+1,\,j+1}-x_{N+j-i,\,j}}\cdot
 e^{x_{N+j-i,\,j+1}-x_{N+j-i+1,\,j+1}}\,,
$$
where $e^{x_{N+j-i-1,\,j}-x_{N+j-i,\,j}}$ with
$e^{x_{N+j-i,\,j+1}-x_{N+j-i+1,\,j+1}}$ come from coefficients in
$E_{i,\,i+1}$,\\ and $e^{x_{N+j-i,\,j+1}-x_{N+j-i-1,\,j}}$ with
$e^{x_{N+j-i+1,\,j+1}-x_{N+j-i,\,j}}$ are from the exponent in
\eqref{GrRightWhittaker}.

\npaa Similarly for the left $(m,N)$-Whittaker vectors one obtains:
 \be
  E_{a+1,\,a}\psi_L^{(m,N)}\,=\,-\frac{1}{\h}\psi_L^{(m,N)}\,,
 \hspace{1.5cm}
  a=m+1,\ldots,N-1\,,
 \ee
since for $E^0_{j+1,\,j}$ being the generators \eqref{GGivRep} with
specialized parameters $\mu_j=\mu_{j+1}=0$:
$$
 E^0_{j+1,\,j}\cdot e^{-\h^{-1}\,e^{x_{N-1,\,j}}}\,
 =\,-\frac{1}{\h}\,e^{-\h^{-1}\,e^{x_{N-1,\,j}}}\,,
\hspace{1.5cm}
 j=1,\ldots,N-1\,,
$$
and
$$
 E^0_{j+1,\,j}\cdot\,\exp\Big\{-\frac{1}{\h}\sum_{i=1}^{r-1}
 e^{x_{j+i,\,j}-x_{j+i-1,\,j}}\,
 +\,e^{x_{j+1+i,\,j+1}-x_{j+i,\,j+1}}\Big\}\,=\,0\,,\quad
 j=1,\,N-1\,,
$$
for any $r=2,\ldots,N-1$ holds due to the "box relations"
(\ref{ToricVariety}.1):
$$
 e^{x_{j+i-1,\,j}-x_{j+i,\,j}}\cdot e^{x_{j+i,\,j+1}-x_{j+i-1,\,j}}\,
 =\,e^{x_{j+1+i,\,j+1}-x_{j+i,\,j}}\cdot
 e^{x_{j+1+i,\,j+1}-x_{j+i,\,j+1}}
 \,,
$$
where $e^{x_{j+i,\,j+1}-x_{j+i-1,\,j}}$ with
$e^{x_{j+1+i,\,j+1}-x_{j+i,\,j}}$ come from coefficients of
$E_{j+1,\,j}$, and $e^{x_{j+i-1,\,j}-x_{j+i,\,j}}$ with
$e^{x_{j+1+i,\,j+1}-x_{j+i,\,j+1}}$ are from the exponent of
\eqref{GrLeftWhittaker}.

\npaa At last to verify the remaining defining relations
$$
 E_{m,\,N}\psi_R^{m,\,N}\,=\,-\h^{-1}\psi_R^{m,\,N}\,,
\hspace{1.5cm}
 E_{m+1,\,1}\psi_L^{m,\,N}\,=\,-\h^{-1}\psi_L^{m,\,N}\,,
$$
for the $(m,\,N)$-Whittaker vectors, we use the explicit form
\eqref{EkN} of the generators $E_{m,\,N}$ and $E_{m+1,\,1}$ and
apply the relations \eqref{PathRelationA}, \eqref{PathRelationAt},
and \eqref{PathRelationB}, \eqref{PathRelationBt}. $\Box$

\section{Proof of Theorem 1.1}

In this part we apply the results of previous Section to derive the
integral formula \eqref{GivIntegralRep} for
$\Psi^{(m,N)}_{\ul}(x,0\ldots,0)$. Namely,  into \eqref{GrWhittaker}
using \eqref{GGivPairing} one readily obtains:
 \be\label{IntStep1}
  \Psi^{(m,N)}_{\ul}(x,0\ldots,0)\,
  =\,e^{-x\frac{m(N-m)}{2}}\bigl\<\psi_L^{(m,N)}\,,\,
  e^{x(E_{11}+\ldots+E_{mm})}\psi_R^{(m,N)}\bigr\>\\
  =\,e^{-x\frac{m(N-m)}{2}}\int_{\CC}\prod_{n=1}^{N-1}\prod_{i=1}^n
  e^{-x_{n,i}}dx_{n,i}\,\,
  \ov{\psi_L^{m,N}(\ux)}e^{x(E_{11}+\ldots+E_{mm})}\psi_R^{(m,N)}(\ux)\,.
 \ee
Then one picks from \eqref{GGivRep}:
$$
 E_{11}+\ldots+E_{mm}\,\,=\,\,\sum_{k=1}^m\mu_k\,\,
 +\,\,\sum_{n=1}^{N-m}\,\,\sum_{i=1}^{\min(m,n)}\!\!
 \frac{\pr}{\pr x_{n,i}}\,\,
 +\,\,\sum_{n=1}^{m-1}\,\,\sum_{i=n+1}^{\min(N-m+n,\,m)}\!\!
 \frac{\pr}{\pr x_{N-m+n,\,i}}\,,
$$
and substituting this together with \eqref{GrLeftWhittaker} and
\eqref{GrRightWhittaker} into \eqref{IntStep1} one finds out the
following:
 \be\label{IntStep2}
  \ldots=\,\frac{1}{\ov{C_{m,N}^L}\,C_{m,N}^R}
  \int\limits_{\CC}\prod_{n=1}^{N-1}\prod_{i=1}^ndx_{ni}\,\,
  e^{\sum\limits_{k=1}^{m-1}\ov{\mu_{N-1-k}}\,x_{N-1-k,\,1}\,
  -\,\sum\limits_{k=m+1}^{N-1}\mu_kx_{kk}}\\
  \times\exp\Big\{\i\sum_{n=1}^{N-1}\bigl(\la_n-\la_{n+1}\bigr)
  \sum_{i=1}^nx_{n,i}\,+\,\i\sum_{k=1}^m\la_{N-m+k}x\,
  -\,\frac{1}{\h}\Big(e^{x_{N-m,\,1}-x}\\
  +\,\sum_{k=1}^m\sum_{i=1}^{N-m-1}e^{x_{i+k-1,\,k}-x_{i+k,\,k}}\,
  +\,\sum_{k=m+1}^{N-1}\Big[e^{x_{N-1,\,k}}\,
  +\,\sum_{i=1}^{N-k-1}e^{x_{i+k-1,\,k}-x_{i+k,\,k}}\Big]\\
  +\,e^{-x_{mm}}\,
  +\,\sum_{k=1}^{N-m}\sum_{i=1}^{m-1}e^{x_{k+i,\,i+1}-x_{k+i-1,\,i}}\\
  +\,\sum_{k=1}^{m-1}\Big[e^{-x_{N-1,\,k}}\,
  +\,\sum_{i=1}^{k-1}e^{x_{N-k+i,\,i+1}-x_{N-k+i-1,\,i}}\Big]
  \Big)\Big\}\,.
 \ee
Next, let us integrate out the variables
$x_{N-k,\,i},\,k=1,\ldots,m-1,\,i=1,\ldots,m-k$ and
$x_{m+k,\,m+j},\,k=1,\ldots,N-m-1,\,i=1,\ldots,k$:
 \be
  \int\limits_{\IR^{\frac{m(m-1)}{2}}}
  \prod_{k=1}^{m-1}\prod_{i=1}^{m-k}dx_{N-k,\,i}\,\,
  \prod_{k=1}^{m-1}\exp\Big\{\i(\la_{N-k}-\la_{N-k+1})
  \sum_{i=1}^{m-k}x_{N-k,\,i}\\
  -\,(\i\la_{N-1-k}+\rho_{N-1-k})x_{N-1-k,\,1}\\
  -\,\frac{1}{\h}\Big(e^{-x_{N-1,\,k}}\,
  +\,\sum_{i=1}^{k-1}e^{x_{N-k+i,\,i+1}-x_{N-k+i-1,\,i}}\Big)
  \Big\}\quad=\quad C^R_{m,N}\,,
 \ee
and
 \be
  \int\limits_{\IR^{\frac{(N-m)(N-m-1)}{2}}}
  \prod_{k=1}^{N-m-1}\prod_{i=1}^kdx_{m+k,\,k+i}\,\,
  \prod_{k=1}^{N-m-1}\exp\Big\{\i(\la_{m+k}-\la_{m+k+1})
  \sum_{i=1}^kx_{m+k,\,m+i}\\
  -\,(\i\la_{m+k}-\rho_{m+k})x_{m+k,\,m+k}\\
  -\,\frac{1}{\h}\Big(e^{x_{N-1,\,m+k}}\,
  +\,\sum_{i=1}^{N-m-k-1}e^{x_{m+i+k-1,\,m+k}-x_{m+i+k,\,m+k}}
  \Big)\Big\}\quad=\quad\ov{C^L_{m,N}}\,.
 \ee
Finally, making cancelations of the normalization constants
$C^R_{m,N}$ and $\ov{C^L_{m,N}}$ in \eqref{IntStep2}, one arrives to
\eqref{GivIntegralRep}, and thus completes the proof of Theorem 1.1.

\section{Proof of Propositions 1.1 and 1.2}

In this part we prove Proposition 1.1. Explicit form
\eqref{GrLaxOperator} of the quantum Lax operator
$\CL(\ux;\,\pr_{\ux})$ readily follows from a simple calculation.
\begin{lem} The adjoint action of $g=g(\ux)\in H^{(m,N)}\subset GL_N$
\eqref{GrGroupElement} in
$\Mat(N,\IR)$ reads as follows:
 \be\label{GrAdjoint}
  g^{-1}E_{11}g\,=\,E_{11}\,-\,\sum_{k=2}^mx_kE_{k,1}\,;\\
  g^{-1}E_{1,k}g\,=\,E_{1,k}\,+\,x_kE_{11}\,
  -\,\sum_{n=2}^mx_n\bigl(E_{n,k}+x_kE_{n,1}\bigr)\,,\quad
  k=2,\ldots,m\,;\\
  g^{-1}E_{k,n}g\,=\,E_{k,n}\,+\,x_nE_{k,1}\,,\quad
  k,n=2,\ldots,m\,;\\
  g^{-1}E_{1,a}g\,=\,e^{x_N-x_1}\Big[E_{1,a}+x_aE_{1,N}\,
  -\,\sum_{k=2}^mx_k\bigl(E_{k,a}+x_aE_{k,N}\bigr)\Big]\,,\\
  g^{-1}E_{k,a}g\,=\,e^{x_N-x_1}\Big[E_{k,a}+x_aE_{k,N}\Big]\,,\quad
  k=2,\ldots,m\,,\,\,a=m+1,\ldots,N-1\,;\\
  g^{-1}E_{a,b}g\,=\,E_{a,b}\,+\,x_bE_{a,N}\,,\quad
  a=m+1,\ldots,N-1\,,\,\,b=a,\ldots,N-1\,.
 \ee
\end{lem}
\emph{Proof:} Expanding the functions
$F_{ij}(\ux)=g(\ux)^{-1}E_{ij}g(\ux)$ by the Taylor formula one
arrives to \eqref{GrAdjoint}. $\Box$

\npa Next, let us introduce the following notation:
$$
 \<Xg\>\,
 =\,\bigl\<\psi_L\,,\,\pi_{\ul}(X)\pi_{\ul}(g)\psi_R\bigr\>\,,
\hspace{1.5cm}
 X\in\CU(\gl_N)\,,\quad g=g(\ux)\in H^{(m,N)}\,.
$$
Thus, using \eqref{GrAdjoint}, and taking into account the property
\eqref{HermiteanProp} with the defining equations
\eqref{LeftWhittEqs}, \eqref{RightWhittEqs}, one can find the
following:
 \be\label{AdjointMatElemnts}
  \<E_{11}g\>\,=\,\<(E_{11}+\ldots+E_{mm})g\>\,-\,\sum_{k=2}^m\<E_{kk}g\>\,
  =\,\pr_{x_1}\<g\>\,;\\
  \<E_{1,k}g\>\,=\,\<E_{1,k}g\>\,+\,x_k\<E_{11}g\>\,
  +\,\sum_{n=2}^mx_n\bigl(x_k\<E_{n,1}g\>-\<E_{n,k}g\>\bigr)\,
  =\,\Big\{-\delta_{k,2}\h^{-1}\\
  +(1-\delta_{k,m})\h^{-1}x_{k+1}\,
  +\,x_k\pr_{x_1}\,+\,\sum_{n=2}^mx_kx_n\pr_{x_n}\Big\}\<g\>\,,\qquad
  k=2,\ldots,m;\\
  \<E_{k,\,i}g\>\,=\,\bigl\{-\delta_{i,\,k+1}\h^{-1}\,
  +\,x_i\pr_{x_k}\bigr\}\<g\>\,,\quad k=2,\ldots,m\,,\,\,
  i=k+1,\ldots,m\,;\\
  \<E_{1,a}g\>\,=\,-(-1)^{\e(m,N)}x_ax_me^{x_N-x_1}\,,\\
  \<E_{m,a}\>\,=\,(-1)^{\e(m,N)}x_ae^{x_N-x_1}\,,\qquad
  a=m+1,\ldots,N-1\,;\\
  \<E_{1N}g\>\,=\,-(-1)^{\e(m,N)}x_me^{x_N-x_1}\,,\qquad
  \<E_{m,N}g\>\,=\,(-1)^{\e(m,N)}e^{x_N-x_1}\,;\\
  \<E_{a,i}g\>\,=\,x_i\pr_{x_a}\,,\quad a=m+1,\ldots,N-1\,,\quad
  i=a,\ldots,N-1\,;\\
  \<E_{NN}g\>\,=\,\Big\{\pr_{x_N}\,
  -\,\sum_{a=m+1}^{N-1}x_a\pr_{x_a}\Big\}\<g\>\,.
 \ee
At last, using \eqref{HermiteanProp}, together with the defining
equations \eqref{LeftWhittEqs}, \eqref{RightWhittEqs} one obtains
the expressions for the remaining matrix elements of the Lax
operator \eqref{GrLaxOperator}.

\npa Finally, let us adopt the following notations:
$$
 C_I\,=\,\sum_{i,j=1\atop i<j}^NE_{ii}E_{jj}\,,
\hspace{1cm}
 C_{II}\,=\,\sum_{i,j=1\atop i<j}^NE_{ji}E_{ij}\,,
\hspace{1cm}
 C_{III}\,=\,\sum_{i=1}^N\rho_iE_{ii}\,;
$$
and therefore (\ref{glCasimirs}.2) reads
$$
 C_2\,=\,C_I\,-\,C_{II}\,-\,C_{III}\,+\,\s_2(\rho)\,.
$$
Then one has
 \be
  C_I\,
  =\,(E_{11}+\ldots+E_{mm})(E_{m+1,\,m+1}+\ldots+E_{NN})\,
  +\,\sum_{j=2}^m(E_{11}+\ldots+E_{mm})E_{jj}\\
  +\,\sum_{i=m+1}^{N-1}E_{ii}(E_{m+1,\,m+1}+\ldots+E_{NN})\,
  -\,\sum_{i,j=2\atop i\leq j}^mE_{ii}E_{jj}\,
  -\,\sum_{i,j=m+1\atop i\leq j}^{N-1}E_{ii}E_{jj}\,;
 \ee
and similarly to \eqref{AdjointMatElemnts} one finds out:
 \be\label{QuadraticCasimir1}
  \<C_Ig\>\,=\,\<C_I'g\>\,
  =\,\Big\{\frac{\pr^2}{\pr x_1\pr x_N}\,
  +\,\sum_{i=m+1}^{N-1}x_i\frac{\pr^2}{\pr x_i\pr x_N}\,
  -\!\!\sum_{i,j=m+1\atop i\leq j}^{N-1}\Big(
  x_i\frac{\pr}{\pr x_i}\Big)\Big(
  x_j\frac{\pr}{\pr x_j}\Big)\Big\}\<g\>\,,
 \ee
with
 \be
  C_I'\,=\,(E_{11}+\ldots+E_{mm})(E_{m+1,\,m+1}+\ldots+E_{NN})\\
  +\,\sum_{i=m+1}^{N-1}E_{ii}(E_{m+1,\,m+1}+\ldots+E_{NN})\,
  -\,\sum_{i,j=m+1\atop i\leq j}^{N-1}E_{ii}E_{jj}\,.
 \ee
Next, let us observe that for $1<k\leq m$ one has
 \be
  \<E_{k,1}E_{1,k}g\>\,
  =\,\Big\{x_k\frac{\pr^2}{\pr x_1\pr x_k}\,
  -(-x_{k-1})^{1-\delta_{k,2}}\frac{\pr}{\pr x_k}\\
  +\,\sum_{n=k}^m\Big[x_n\frac{\pr}{\pr x_n}\,
  +\,\Big(x_k\frac{\pr}{\pr x_k}\Big)
  \Big(x_n\frac{\pr}{\pr x_n}\Big)\Big]\Big\}\<g\>\,;
 \ee
hence, for
$$
 \<C_{II}g\>\,
 =\,
 \sum_{k=2}^m\<E_{k,1}E_{1,k}g\>\,
 +\,\sum_{i=m+1}^{N-1}\<E_{i+1,\,i}E_{i,\,i+1}\>\,
 +\,\<E_{m+1,\,1}E_{1,\,m+1}g\>\,,
$$
one obtains
 \be\label{QuadraticCasimir2}
  \sum_{k=2}^m\<E_{k,1}E_{1,k}g\>\,
  =\,\Big\{\h^{-1}\Big((\delta_{m,1}-1)\frac{\pr}{\pr x_2}\,
  +\sum_{k=2}^{m-1}x_k\frac{\pr}{\pr x_{k+1}}\Big)\\
  +\,\sum_{k=2}^m\Big[(k-1)x_k\frac{\pr}{\pr x_k}\,
  +\,\frac{\pr^2}{\pr x_k^2}\Big]\,
  +\,\sum_{i,j=1\atop i<j}^m
  \Big(x_i^{1-\delta_{i,1}}\frac{\pr}{\pr x_i}\Big)
  \Big(x_j\frac{\pr}{\pr x_j}\Big)\Big\}\<g\>\,,\\
  \sum_{i=m+1}^{N-1}\<E_{i+1,\,i}E_{i,\,i+1}\>\,
  =\,\Big\{\sum_{i=m+1}^{N-2}x_{i+1}\frac{\pr}{\pr x_i}\,
  +\,(1-\delta_{m,\,N-1})\h^{-1}\frac{\pr}{\pr x_{N-1}}
  \Big\}\<g\>\,,\\
  \<E_{m+1,\,1}E_{1,\,m+1}g\>\,
  =\,(-1)^{\delta_{m,\,N-1}+\e(m,N)}\h^{-2}\,
  (x_m)^{1-\delta_{m,1}}
  (x_{m+1})^{1-\delta_{m,\,N-1}}\,e^{x_N-x_1}\,\<g\>\,.
 \ee

At last one derives
 \be\label{QuadraticCasimir3}
  \<C_{III}g\>\,
  =\,\Big\{\rho_1\Big(\frac{\pr}{\pr x_1}-\frac{\pr}{\pr x_N}\Big)\,
  +\,\sum_{i=m+1}^{N-1}(\rho_1+\rho_i)x_i\frac{\pr}{\pr x_i}
  \Big\}\<g\>\,,
 \ee
and collecting \eqref{QuadraticCasimir1}, \eqref{QuadraticCasimir2}
and \eqref{QuadraticCasimir3} one arrives at
(\ref{GrHamiltonians}.2).

\vspace{1cm}

\noindent
{\sc Institute for Theoretical and Experimental Physics,\\
Bol. Cheremushkinskaya 25, Moscow 117218,\\
\emph{E-mail address:}\quad \verb"Sergey.Oblezin@itep.ru"}

\end{document}